\documentclass{amsart}
\usepackage{amssymb}

\newtheorem{thm}{Theorem}[section]
\newtheorem{cor}[thm]{Corollary}
\newtheorem{prop}[thm]{Proposition}

\theoremstyle{definition}
\newtheorem{defi}[thm]{Definition}
\newtheorem{ex}[thm]{Example}
\newtheorem{rem}[thm]{Remark}

\numberwithin{equation}{section}

\newcommand{\C}{\mathbb C}   
\newcommand{\R}{\mathbb R}
\newcommand{\Z}{\mathbb Z}

\begin{document}

\title[Rigidity of Variations of Hodge Structure]{Coupled Contact Systems and Rigidity of Maximal Dimensional Variations of Hodge Structure}

\author{Rich\'ard Mayer}
\address{Department of Mathematics, University of Massachusetts, Amherst, MA 01003}
\email{mayer@math.umass.edu}
\subjclass{Primary 14C30}

\begin{abstract}
In this article we  prove that locally Griffiths' horizontal distribution on the period domain is given by a generalized version of the familiar contact differential system. As a consequence of this description we obtain strong local rigidity properties of maximal dimensional variations of Hodge structure. For example, we prove that if the weight is odd (greater than one) then there is a unique germ of maximal dimensional variation of Hodge structure through every point of the period domain. Similar results hold if the weight is even with the exception of one case.
\end{abstract}

\maketitle

\section{Introduction}

Variations of Hodge structure are integral manifolds of Griffiths' horizontal distribution on the period domain \(D\) of a fixed integral lattice which are stable under the action of a discrete subgroup of the isometry group of \(D\). This distribution is defined in terms of the filtrations \(\{F^p\}\) by 
\[ 
\frac{\partial F^p}{\partial z_i}\subset F^{p-1}.
\]
The horizontal distribution is not completely integrable in general if the weight is greater than one, so integral manifolds must have lower dimension than the distribution. There are two fundamental questions for not completely integrable differential systems: 
\begin{itemize}
\item What is the maximal dimension an integral manifold can have? \\
\item What is the structure of the integral manifolds that attain the dimension bound? 
\end{itemize}
We shall examine these two questions for the horizontal system in detail. The tool that makes this possible is the fact that the horizontal differential system is locally given by coupled matrix valued contact systems. (See Section 4.) We derive two consequences of this description. The first one is a new proof of the main result of \cite{C-K-T} which says that there are explicit quadratic functions \(q^{even}_1, q^{even}_2, q^{even}_3\) (weight even) and \(q^{odd}_1, q^{odd}_2\) (weight odd) of the Hodge numbers \(h^{i,j}\) which give sharp upper bounds for the dimension of a variation of Hodge structure.
(Theorem~\ref{dimboundthm}). In fact, we cover a case which is missing from \cite{C-K-T}, this will be explained in Section 4. The second consequence is the main result of this work which generalizes Carlson's rigidity theorem for weight two variations (cf. \cite{C1}) to arbitrary weights. The precise statement is as follows:
\begin{thm}\label{mainthm}
Let \(D\) denote the period domain and let \(w\) be the weight. 
\begin{enumerate} 
\item Assume that one of the following holds:
\begin{enumerate}
\item \(w=2k+1\), \(h^{k,k+1}>2\) and all the other Hodge numbers are greater than one
\item \(w=2k\), all the Hodge numbers are greater than one and the maximum dimension is \(q^{even}_1\) 
\end{enumerate}
then there is a unique maximal dimensional germ of variations of Hodge structure through each point of \(D\), i.e., there is a unique maximal dimensional foliation the leaves of which are the maximal dimensional variations of Hodge structure.
\item Let \(S_1, S_2\subset D\) be two maximal dimensional variations of Hodge structure. Assume that one of the following holds:
\begin{enumerate}
\item \(w=2k\), \(dim(S_1)=dim(S_2)=q^{even}_2\), \(h^{k+1}>2\), the other Hodge numbers are greater than one and \(h^k\geq 4\) is even
\item \(w=2k\), \(dim(S_1)=dim(S_2)=q^{even}_3\) and \(h^{k+1}>2\) 
\end{enumerate}
then there is an element \(g\in Aut(D)\) such that \(g\cdot S_1=S_2\) holds in a neighborhood.
\item Let \(w=2k\) and let \(E\) denote a maximal dimensional integral element. If \(dim(E)=q^{even}_2\) and \(h^k\) is odd then there is an infinite dimensional family of germs of maximal dimensional integral manifolds tangent to \(E\), i.e. flexibility holds.
\end{enumerate}
\end{thm} 

In the remaining special cases when the Hodge numbers are small (i.e., they do not satisfy the requirements of the theorem) rigidity does not hold. In fact, in these cases we have the behavior of the classical contact system explained in Section~\ref{contactsection}. 

The outline of the paper is as follows. In Section 2 we give the necessary definitions of exterior differential systems and treat the case of the classical contact system. In Section 3 we introduce local coordinate systems on the period domain. Using these coordinates we show in Section 4 that the horizontal system is locally equivalent to a system given by analogs of the classical contact system. Analyzing this new system leads to the proof of the main theorem in Section 5. 

The author would like to thank his advisor James A. Carlson for suggesting this problem to him and for the support throughout his graduate studies.

\section{Exterior Differential Systems}

The main reference for this section is \cite{B}. Let \(X\) be a manifold, let \(T^*\) denote the cotangent sheaf of \(X\) and \(\bigwedge ^* T^*\) the associated deRham algebra. We will be concerned with complex manifolds and the holomorphic cotangent sheaf, but most of the results remain true even if we consider \(C^\infty \) real manifolds. 
\subsection{Integral Elements and Manifolds}
\begin{defi}
A differential system is defined by an ideal \(\mathcal{I} \subset \bigwedge ^* T^*\) which is closed under exterior differentiation, i.e., \(d\mathcal{I}\subset\mathcal{I}\). An integral manifold of \(\mathcal{I}\) is a holomorphic mapping \(i\):~ \(S \longrightarrow X\) such that \(i^*\omega = 0\) for each germ \(\omega\) of \(\mathcal{I}\).
\end{defi} 

Often we will take \(i\) to be the inclusion and we will talk about the {\it germ of an integral manifolds} by which we will mean an equivalence class of integral manifolds where two integral manifolds are equivalent if their images agree in a neighborhood of a point. 

Let \(i\):~ \(S \longrightarrow X\) be an integral manifold of \(\mathcal{I}\). If \(s \in S \) and \(E = T_s S \subset T_s M\) is the tangent space to \(S\) at \(s\) then for each \(\omega \in \mathcal{I}\) \(\omega_E = 0\) where \(\omega_E\) denotes the restriction of \(\omega\) to \(E\). It is clear that the vanishing of \(i^*\omega\) at each point depends only on the tangent space of \(S\). This leads to the following definition.

\begin{defi}
A linear subspace \(E \subset M\) is an integral element of \(\mathcal{I}\) if \(\omega_E = 0\) for each \(\omega \in \mathcal{I}\). 
\end{defi}

Note that in general an integral element is not necessarily tangent to an integral manifold, see \cite{B} for examples. 

Suppose now that \(\mathcal{I}\)  is generated algebraically by the
1-forms \(\{\omega^1, \ldots, \omega^{n-k}\}\) where \(n\) is the
dimension of \(M\) and let \({\mathcal{D}}_{v} \subset T_{v} M\) denote
the linear space of tangent vectors that are annihilated by \(\mathcal{I}\). If \({\mathcal{D}}_{v}\) has the same dimension for each \(v\in M\) then these subspaces define a {\it distribution} in the tangent bundle. The dimension of the distribution \(\mathcal{D}\) is the dimension of \({\mathcal{D}}_{v}\). From now on we will consider differential systems that are defined by 1-forms of this type and we will use the words differential system and distribution interchangeably.

 The condition that \(\mathcal{I}\) is closed means that \(d\omega^i\) is
 given as an algebraic combination of \(\{\omega^1, \ldots,
 \omega^{n-k}\}\) for each \(i\). This condition \((F)\) is called the
 {\it Frobenius condition}. It is not hard to see (cf. \cite{W}
 Proposition 2.30) that the Frobenius condition for differential forms is equivalent to the following for a distribution \(\mathcal{D}\): for any two holomorphic vector fields \(X\) and \(Y\) lying in \(\mathcal{D}\), the Lie bracket [\(X,Y\)] also lies in \(\mathcal{D}\). The fundamental integrability result for such systems is the following theorem of Frobenius.

\begin{thm}[Frobenius]\label{frobenius}
Let \(\mathcal{I}\) be a differential ideal generated algebraically by the linearly independent 1-forms \(\{\omega^1, \ldots, \omega^{n-k}\}\). There is a local coordinate system \((y^1,\ldots ,y^n)\) at each point of \(M\) such that \(\mathcal{I}\) is generated by \(\{dy^{k+1}, \ldots, dy^{n}\}\) if and only if the condition \((F)\) holds.
\end{thm}

If such coordinate system exists then clearly there is a \(k\)-dimensional integral manifold through each point of \(M\). For this reason such differential systems are called {\it completely integrable}. In this work we will be interested in differential systems that are not completely integrable. If this is the case then the maximal possible dimension of an integral manifold must be strictly less than the dimension of the distribution. It is an interesting question to determine this dimension. 
\subsection{Contact Differential Systems}\label{contactsection}

Let us start with a differential system defined by a single 1-form.  

\begin{defi}
Let \(M=\C^{2n+1}\) with coordinates \((x_1,\ldots,x_n,y_1,\ldots,y_n,z)\).
The differential system whose differential ideal is generated by the 1-form
\[ \omega = dz - \sum_{i=1}^{n}x_idy_i \]
is called the contact system.
\end{defi}
\begin{rem}
Sometimes we will refer to this system as the classical contact system
in order to distinguish it from matrix valued analogs that will be
discussed later.
\end{rem}
Let us observe that the distribution defined by this system has dimension \(2n\). However, the system is not completely integrable (as we will see in a moment) so any integral manifold must have dimension less than \(2n\). An easy computation shows that 
\[\omega \wedge (dw)^{\wedge n} = \pm dz\wedge dx_1\wedge \ldots \wedge dx_n \wedge dy_1 \wedge \ldots \wedge dy_n \] 
so \(d\omega \wedge \omega\) is clearly nonzero which implies that the Frobenius condition can not hold for the contact system. In light of this, our next task is going to be to determine the maximal dimension an integral manifold can have. 

In what follows the notation \(f(x_I,y_J)\) means that the function \(f\) depends on the variables \(x_i, i \in I\) and \(y_j, j \in J\).

\begin{prop}\label{contactmax}
Let \(S \subset {\C}^{2n+1}\) be a \(k\)-dimensional integral manifold
of the contact system going through the point \(s\). Then there exist
two disjoint sets of indices \(I,J\subset \{1,\ldots,n\}\) of total
length \(k\) such that in a neighborhood of \(s\), \(S\) is
parameterized by a set of holomorphic functions
\((g_1(x_I,y_J),\ldots,g_{2n+1}(x_I,y_J))\) where the variables
\(x_{i}\) and \(y_{j}\) \(({i\in I},{j \in J})\) are independent on \(S\). 
\end{prop}
\begin{rem}
By definition the variables \(x_{i}\), \(y_{j}\) \(({i\in I},{j \in
J})\) are independent on \(S\) if the differential 1-forms
\(dx_{i}\) and \(dy_{j}\) \(({i\in I},{j \in J})\) are linearly independent at each point of a
neighborhood of
\(S\). This also means that independent variables form part of a
coordinate system in a neighborhood. 
\end{rem}
\begin{proof}
In a neighborhood \(U\) of \(s\) the manifold \(S\) is given by a set
of holomorphic functions \((g_1,\ldots,g_{2n+1})\). Let \(L\subset \{x_i, y_j,z\}\) denote
the set of {\em independent} coordinates that the functions \(g_i\) depend on,
(\(|L|=k\)). Let \(I\subset L\), \(J\subset L\) be the sets consisting of
variables \(x_{i}\), \(y_{j}\) respectively. Since \(S\) is an integral manifold, \(\omega\) and \(d\omega\) must vanish on \(S\). This implies that the coordinate \(z\) can be expressed as a function of the other variables so we can assume that the functions \(g_i\) do not depend on \(z\), i.e., \(z\notin L\).

Suppose that there is a pair of coordinates \((x_i,y_i)\), \(x_i,y_i\in L\) such that \(dx_i \wedge dy_i \) does not vanish on \(U\). Since 
\[
0 = d\omega = \sum_{l=1}^{n}dx_l \wedge dy_l \;\;\;\;\;\;\mathrm{on}\;\;\;S
\]
there must be another index \(j\neq i\)  such that \(dx_j \wedge dy_j
\) does not vanish on \(U\) and the variables \(x_j,y_j\) are not in \(L\) since they depend on \(x_i,y_i\). Now 
\[0 \neq dx_j\wedge dy_j = (\frac{\partial x_j}{\partial x_i} \frac{\partial y_j}{\partial y_i} - \frac{\partial y_j}{\partial x_i} \frac{\partial x_j}{\partial y_i})dx_i\wedge dy_i + \eta\]
where \(\eta\) is a 2-form that does not contain \(dx_i\wedge dy_i\). This means that one of the coefficient functions, say \(\frac{\partial x_j}{\partial x_i}\), must be nonvanishing so we can introduce the coordinate change \(x_i \longleftrightarrow x_j\) in a neighborhood of \(s\), by which we mean that we replace the coordinate \(x_i\) by \(x_j\) and leave the other coordinates intact. This means that we replaced the independent pair of variables \((x_i,y_i)\) by another independent pair \((x_j,y_i)\) thereby decreasing \(|I\cap J|\). Continuing in the same manner, after at most \(k\) steps we arrive at two {\em disjoint} sets of indices \(I,J\subset \{1,\ldots,n\}\) that satisfy the requirements of the proposition.  
\end{proof}
\begin{rem}\label{fixj}
It is clear from the proof that if we choose the
set \(J\) to be maximal in the sense that no other independent set of variables that \(S\) depends on contains more \(y_j\)'s then we
can find the set \(I\) with \(|I\cap J|=0\) without changing the set
\(J\). This simple remark will be used quite frequently later,
sometimes without explicit reference to it. 
\end{rem}

\begin{cor}
Let \(S \subset {\C}^{2n+1}\) be an integral manifold of the contact system. Then \(dim(S) \leq n\).
\end{cor}

\begin{proof}
By Proposition~\ref{contactmax}, \(S\) is locally defined by functions \(g_i\) depending on at most \(n\) independent variables. This clearly implies the claim.
\end{proof}

The following
result of \cite{A} (Appendix 4) gives a complete characterization of maximal dimensional
integral manifolds of the contact system.

\begin{thm}\label{arnold}
For any partition \(I+J\) of the set of indices \(\{1,\ldots,n\}\) into two disjoint subsets and for any function \(f(x_I,y_J)\) of \(n\) variables \(x_i,y_j\) \((i\in I, j\in J)\), the formulas 
\[y_i = \frac{\partial f}{\partial x_i} \; (i\in I),\; \; \;  x_j = -\frac{\partial f}{\partial y_j}\; (j\in J),\; \; \; z = f - \sum_{i\in I}x_i\frac{\partial f}{\partial x_i} \]
define an \(n\)-dimensional integral manifold of the contact system.

Conversely, every \(n\)-dimensional integral manifold is defined in a neighborhood of every point by these formulas for a choice of the subset \(I\) and for some generating function \(f\). 
\end{thm}

\begin{proof}
A simple calculation shows that the manifold defined by the formulas in the first part of the Theorem is indeed an integral manifold of the contact system. 

Assume now that \(S\) is an \(n\)-dimensional integral manifold. By Proposition~\ref{contactmax} we can find \(n\) {\em independent} variables \(x_I,y_J\) such that \(S\) is defined by functions of these variables. Let 
\[f = z + \sum_{i\in I} x_iy_i\]
and notice that 
\[df = \sum_{i\in I}y_idx_i - \sum_{j\in J}x_jdy_j \]
since \(0 = \omega = dz + \sum x_ldy_l\) on \(S\). Also, the differential forms \(dx_i,\, dy_j\) are independent so it follows from the defining equation of \(f\) that 
\[y_i = \frac{\partial f}{\partial x_i} \; (i\in I),\; \; \;  x_j = -\frac{\partial f}{\partial y_j}\; (j\in J)\]
which finishes the proof. 
\end{proof}
\begin{rem}\label{contactflexible}
Let us mention here that if we fix an \(n\)-dimensional integral element at the point \(s\), then there exists an integral manifold that is tangent to the given integral element (see Theorem~\ref{exptrick}) but this integral manifold is by no means unique. This follows from the fact that we can choose \(f\) in the previous theorem to be an arbitrary holomorphic function and fixing an integral element specifies only the linear part of \(f\). This means that to any integral element of dimension \(n\) there is an infinite dimensional family of germs of integral manifolds that are tangent to the integral element at the point \(s\).  
\end{rem}
%%%%%%%%%%%%%%
\section{Variations of Hodge Structure}
In this section we will consider variations of Hodge structure as integral manifolds of the horizontal distribution. First let us recall the necessary definitions.
\subsection{Hodge Structures and Classifying Spaces}

For a more detailed discussion, proofs of the basic results, as well as geometric motivation for some of the following definitions, see \cite{G1} and \cite{G2}. 

Let us fix a finite dimensional real vector space \(H_{\R}\) and a lattice \(H_{\Z}\subset H_{\R}\) together with an integer \(w\), (\(w\) will be referred to as the weight). Suppose that we are given a non-degenerate bilinear form \(Q\) on \(H_{\R}\) which satisfies \(Q(x,y) = (-1)^w Q(y,x)\) for all \(x,y\in H_{\R}\) and takes rational values on \(H_{\Z}\). Let \(H_{\C} = H_{\R} \otimes_{\R} \C\) denote the complexification of \(H_{\R}\). 

\begin{defi}
A Hodge structure of weight \(w\) on \(H_{\R}\) is a decomposition  
\[H_{\C} = \bigoplus_{p+q = w}H^{p,q}\] 
such that \(H^{p,q}\) and \(H^{q,p}\) are complex conjugate to each other with respect to \(H_{\R}\). The integers \(h^{p,q} = dim\,H^{p,q}\) are the Hodge numbers.
\end{defi}

Let \(h(x,y) = (-1)^{w(w-1)/2}Q(x,\overline{y})\) and let \(C\) denote the Weil operator that acts on \(H^{p,q}\) by multiplication by \(i^{p-q}\).

\begin{defi}
The Hodge structure is weakly polarized if 
\[Q(H^{p,q},H^{r,s})=0 \; \; \; \;\rm{for} \;\;\;\; (r,s) \neq (q,p).\]
It is strongly polarized if, in addition, the form \(h_C(x,y) = h(Cx,y)\) is positive hermitian.
\end{defi}
%%%%filtrations 
To each Hodge structure of weight \(w\) we can assign the {\it Hodge filtration}
\[H_{\C} \supset \ldots \supset F^{p-1} \supset F^{p} \supset F^{p+1}\supset \ldots \supset 0\]
where 
\[F^{p} = \bigoplus_{i\geq p}H^{i,w-i}.\]
This filtration has the property
\[H_{\C} = F^{p} \oplus \bar{F}^{w-p+1} \; \; \; \rm{for}\; \;\rm{each} \; \;p\]
where the bar denotes complex conjugation.
If the Hodge structure is strongly polarized then the associated Hodge filtration satisfies \(Q(F^p,F^{w-p+1})=0\) for all \(p\), and the form \(h_C(x,y)\) is positive hermitian. It is not hard too see that a Hodge filtration with the listed properties determines a strongly polarized Hodge structure so these two notions are in fact equivalent. 

Let \(\check{D}\) denote the set of weakly polarized Hodge structures (or Hodge filtrations) with fixed Hodge numbers \(h^{p,q}\), \((p+q=w)\) and let \(D \subset \check{D}\) be the subset of strongly polarized Hodge structures. There is a natural complex structure on the set \(\check{D}\) which can be described as follows. Since a Hodge structure is determined by its associated filtration, we can view \(\check{D}\) as a subset of a product of Grassmannian manifolds. The conditions \(F^{p} \supset F^{p+1}\) and \(Q(F^p,F^{w-p+1})=0\) are algebraic so \(\check{D}\) is a subvariety of a product of Grassmannians. As such it has the structure of a complex projective variety. It can be checked that the special orthogonal group \(G_{\C} = SO(Q,\C)\) acts transitively on \(\check{D}\), in particular \(\check{D}\) is a complex manifold. The subset \(D\subset \check{D}\) is open in the complex topology and it is homogeneous for the corresponding real group \(G_{\R}\). The classifying space \(D\) is sometimes referred to as the {\it period domain} of Hodge structures.

\subsection{Canonical Coordinates on the Period Domain}

In this section we will describe two sets of local coordinate systems on the period domain. These coordinate systems will be used to investigate the local properties of differential systems on the period domain. The two coordinate systems are equivalent to each other and the reason they will both be used is that certain results can be described more conveniently in one of these systems.
\subsubsection{Lie algebra coordinates}
Let \(\mathfrak{g}_{\C}\) denote the Lie algebra of \(G_{\C}\) and let us fix a reference Hodge structure \(H\in D\). \(H\) defines a Hodge structure of weight 0 on \(\mathfrak{g}_{\C}\) (cf. \cite{S}) where 
\[\mathfrak{g}^{p,-p} = \{\phi \in \mathfrak{g}_{\C}\,|\, \phi(H^{r,s})\subset H^{r+p,s-p}\; \rm{for \;all} \;(r,s)\}.\]
\(\mathfrak{g}_{\C}\) can be decomposed as 
\[\mathfrak{g}_{\C} = \mathfrak{g}^{-}\oplus \mathfrak{g}^{0}\oplus \mathfrak{g}^{+}\]
where the subalgebras \(\mathfrak{g}^{-}\), \(\mathfrak{g}^{0}\) and \(\mathfrak{g}^{+}\) are defined as 
\[\mathfrak{g}^{-} = \bigoplus_{p<0}\mathfrak{g}^{p,-p},\;\;\;\;
  \mathfrak{g}^{0} = \mathfrak{g}^{0,0},\;\;\;\;
  \mathfrak{g}^{+} = \bigoplus_{p>0}\mathfrak{g}^{p,-p}.\]
Recall that \(\check{D} \cong G_{\C}/B\) as a homogeneous space where
\[B = \{g\in G_{\C} \,|\, g(F^{p}) \subset F^{p}\}\]
is the isotropy subgroup of the reference Hodge structure. The Lie algebra of \(B\) is 
\[\mathfrak{b} = \mathfrak{g}^{0}\oplus \mathfrak{g}^{+}\]
so the complement \(\mathfrak{g}^{-}\) can be identified with the holomorphic tangent space of \(\check{D}\). Note that the Lie bracket is compatible with the Hodge decomposition, i.e.,
\begin{equation}\label{liecompat}
[\mathfrak{g}^{p,q},\mathfrak{g}^{r,s}]\subset \mathfrak{g}^{p+r,q+s}.
\end{equation}

\begin{defi}\label{algebracoord}
The local Lie algebra coordinates in a neighborhood of the reference Hodge structure \(H\in D\) are given by the map 
\[\Phi:\mathfrak{g}^{-}\longrightarrow D,\;\;\;\; N\longmapsto (exp(N))\cdot H\]
\end{defi}

Note that the image will lie in \(D\) if the norm of \(N\) is small since \(D\subset \check{D}\) is an open subset. Also, the differential of the map is the identity so \(\Phi\) defines local coordinates in a neighborhood. The actual local coordinates will be chosen by specifying a convenient basis in \(\mathfrak{g}^{-}\). This is our next task (cf. \cite{C-K-T}).
\subsubsection{Hodge Frames and Block Decomposition}

\begin{defi}
A Hodge frame for \(H\) is a set of bases 
\(B^{p,q} = \{B^{p,q}_j \,|\, j=1,\ldots,h^{p,q}\}\)
such that 
\begin{enumerate}
\item[(i)] \(B^{p,q}\) is an \(h_C\)-unitary basis of \(H^{p,q}\)
\item[(ii)] \(B^{p,q} = \overline{B^{q,p}}.\)
\end{enumerate}
\end{defi}

It is clear that the matrix of the bilinear form \(Q\) relative to a Hodge frame has a block decomposition such that the only nonzero blocks are on the antidiagonal and these blocks are identity matrices up to a sign. Let us denote by \(M[i,j]\) the matrix whose only nonzero block is the matrix \(M\) of size \(h^{w-i,i}\times h^{w-j,j}\) in the \((i,j)\) position. Then the matrix of the polarization \(Q\) can be written as 
\[Q=\sum_{k=0}^{w} (-1)^kI[k,w-k]\]
where \(I\) stands for the identity matrix. 

Similarly, the matrix of an endomorphism \(X\) of \(H_{\C}\) relative to a Hodge frame also has a block decomposition \(\sum X_{i,j}[i,j]\) where the matrix \(X_{i,j}\) represents an endomorphism from \(H^{w-j,j}\) to \(H^{w-i,i}\). Now let us consider the block decomposition of a general element \(X\in \mathfrak{g}^{-}\). Clearly the matrix of \(X\) is strictly lower triangular and because it is and element of the orthogonal Lie algebra it satisfies the relation \(X^tQ+QX=0\). In terms of the \(X_{i,j}\) we have
\begin{equation}\label{ortcond} 
(-1)^{w-i}(X_{i,j})^t+(-1)^{w-j}X_{w-j,w-i}=0\;\;\;\;i,j\in \{0,\ldots ,w\}.
\end{equation}
\begin{ex}
To illustrate the definitions above let us consider the weight two case. The matrix of \(Q\) will be 
\[Q=
\begin{pmatrix}
	0&0&I_{h^{2,0}}\\
	0&-I_{h^{1,1}}&0\\
	I_{h^{2,0}}&0&0
\end{pmatrix}
\]
and a general element \(X\in \mathfrak{g}^{-}\) can be written as 
\[X=
\begin{pmatrix}
	0&0&0\\
	X_{1,0}&0&0\\
	X_{2,0}&X_{2,1}&0
\end{pmatrix}
\]
where \((X_{1,0})^t=X_{2,1}\) and \((X_{2,0})^t=-X_{2,0}\) because of Equation~\ref{ortcond}. In this case \(dim(D) = h^{2,0}h^{1,1} + \frac{1}{2}h^{2,0}(h^{2,0}-1)\) since this is the number of coordinates in the Lie algebra coordinate system given by the entries of the matrices \(X_{i,j}\).
\end{ex}

\subsubsection{Lie Group Coordinates}\label{groupcord} 
Now it is easy to describe the other coordinate system. Let \(G^-\) denote the unipotent Lie group which corresponds to the Lie algebra \(\mathfrak{g}^{-}\). 
\begin{defi}\label{groupcoord}  
The local Lie group coordinates in a neighborhood of the reference Hodge structure \(H\in D\) are given by the map 
\[\Psi:G^{-}\longrightarrow D,\;\;\;\; g\longmapsto g\cdot H\]
\end{defi}

It follows from the discussion of the Lie algebra coordinates that this is a coordinate system in a neighborhood of \(H\). If we fix a Hodge frame of \(H\) then an element \(Y\in G^-\) has a block decomposition \(\sum Y_{i,j}[i,j]\) which is induced by the corresponding Lie algebra element \(log(Y)\). 

\begin{rem}
Let us note that the blocks \(Y_{i,j}\) will satisfy the condition
corresponding to Equation~\ref{ortcond} and we will analyze this
later. 
\end{rem}
\begin{rem}
To visualize the matrices in the Lie group coordinate system see Example~\ref{visex}.
\end{rem}
\subsection{Variations of Hodge Structure}

\subsubsection{Horizontal Distribution}

In order to define a distribution we have to specify a fixed dimensional linear subspace of the holomorphic tangent space at each point of \(\check{D}\). The holomorphic tangent space was identified by \(\mathfrak{g}^{-}\) so we can proceed as follows.
\begin{defi}
The horizontal distribution on \(\check{D}\) is given by the vector space \(\mathfrak{g}^{-1,1}\subset \mathfrak{g}^{-}\) at each point of \(\check{D}\).
\end{defi}

This distribution is holomorphic and homogeneous as can be seen by elementary Lie group theory (cf. \cite{S}). Integral manifolds of the horizontal distribution are sometimes referred to as {\it horizontal maps}. We have the following 
\begin{prop}\label{nonintegrable}
The horizontal differential system is not completely integrable provided that \(\mathfrak{g}^{-2,2}\neq 0\).
\end{prop}
\begin{proof}
Since \([\mathfrak{g}^{-1,1},\mathfrak{g}^{-1,1}]=\mathfrak{g}^{-2,2}\), the Proposition follows from Theorem~\ref{frobenius}.
\end{proof}
\begin{defi}
Let \(\Gamma\) denote a properly discontinuous group of automorphisms of \(D\) and let \(M\) be a complex manifold. A holomorphic map \(\Phi : M \longrightarrow \Gamma~\backslash~ D\) is a variation of Hodge structure if the map \(\Phi\) is locally liftable to \(D\) and the local liftings are horizontal (i.e., integral manifolds of the horizontal differential system).
\end{defi}

In this work we will be interested in the local properties of variations of Hodge structure.

\begin{defi}
A germ of a variation of Hodge structure is an equivalence class of integral manifolds \(\Phi :U \longrightarrow D\) of the horizontal system, where two integral manifolds are equivalent if they have the same image in a neighborhood of a point of \(D\).
\end{defi}

\subsubsection{Local Description of the Horizontal Differential System}

Let us consider the description of the horizontal differential system in terms of the Lie algebra and Lie group coordinate systems. 

At the reference Hodge structure the horizontal system is defined by
the differential 1-form entries of the matrices that are at least two
steps below the diagonal
\begin{equation}\label{systemorigin}
dX_{i,j}\;\;\;\;\rm{for}\;\;\;\; i\geq j+2.
\end{equation}
This is clear since at the reference point the horizontal distribution is defined by the subspace \(\mathfrak{g}^{-1,1}\) and the differential forms in Equation~\ref{systemorigin} give the dual of this subspace. Let us now determine the differential forms that define the horizontal differential system in a neighborhood of the reference structure. Since the horizontal system is homogeneous it is given by the left invariant extensions of the differential forms in Equation~\ref{systemorigin}. Let \(\Omega\) denote the Maurer-Cartan matrix, which can be written as \(\Omega = exp(-X)d\,exp(X)\) or \(\Omega = Y^{-1}d\,Y\) in the Lie algebra and Lie group coordinate systems, respectively. Then, in light of the above discussion, we have the following.
\begin{prop}\label{maurer}
The horizontal differential system in a neighborhood of the reference Hodge structure is given by the differential form entries of the matrices
\begin{equation}\label{maurerequation}
\Omega_{i,j}\;\;\;\;\rm{for}\;\;\;\; i\geq j+2.
\end{equation}
\end{prop}

\subsubsection{Integral Elements of the Horizontal Differential System}

Now that we have a local description of the horizontal system we can examine the question: for which integral elements can we find an integral manifold tangent to it? First let us give a characterization of the integral elements.
\begin{prop}\label{commutative}
Let \(E\subset \mathfrak{g}^{-1,1}\) be a linear subspace. Then \(E\) is an integral element of the horizontal differential system if and only if it is an abelian subalgebra.
\end{prop}
\begin{proof}
See \cite{C-K-T} Proposition 2.2 or Section~\ref{proofcom} for another proof.
\end{proof}

Using this result we find the following theorem.
\begin{thm}\label{exptrick}
Let \(E\subset \mathfrak{g}^{-1,1}\) be an integral element at \(v\in D\). Then there is a germ of an integral manifold through \(v\) which is tangent to \(E\).
\end{thm}
\begin{proof}
Let \((e_1,\ldots,e_k)\) be a basis of \(E\). According to Proposition~\ref{commutative} the Lie brackets 
\begin{equation}\label{bracketcomm}
[e_i,e_j]=0\;\;\;\; \rm{for} \;\;\;\; i,j\in (1,\ldots,k)
\end{equation}
 and since \((e_1,\ldots,e_k)\subset \mathfrak{g}^{-1,1}\) the basis elements \(e_i\) are represented by matrices with all blocks zero except the ones one step below the diagonal. Our aim is to show that the following map \(\varphi\) is an integral manifold in a neighborhood \(U\) of the origin of \(\C^{k}\). 
\begin{equation}
\varphi : \C^{k} \longrightarrow D,\;\;\;\;
(x_1,\ldots,x_k)\longmapsto \Phi (x_1\cdot e_1,\ldots,x_k\cdot e_k)
\end{equation}
where the map \(\Phi\) is the Lie algebra coordinate system given in Definition~\ref{algebracoord}. Note that the map \(\varphi\) is well defined and injective in a neighborhood of the origin and its tangent space at the origin is \(E\). What remains to be shown is that \(\varphi\) is horizontal. By Proposition~\ref{maurer} it needs to be checked that the entries \(\Omega_{i,j}\;\;\rm{for}\;\; i\geq j+2\) of the Maurer-Cartan matrix vanish along the image of \(\varphi\). Now 
\[\Omega |\varphi (x_1,\ldots,x_n) = exp(-\sum x_i\cdot e_i)\;d\,exp(\sum x_i\cdot e_i)=\]
\[=\prod exp(-x_i\cdot e_i)\;d\,\prod exp(x_i\cdot e_i)=e_1dx_1+\ldots +e_kdx_k\]
where the second equality holds because of Equation~\ref{bracketcomm}. Since the blocks \((e_l)_{i,j} = 0\) for \(i\geq j+2\) the proposition follows.
\end{proof}

\section{Coupled Contact Systems and Dimension Bounds}

In the previous section we gave a local description of the horizontal differential system and we saw that it is given by the 1-form entries in those blocks of the Maurer-Cartan matrix that are at least two steps below the diagonal (Equation~\ref{maurerequation}). Now we will see that this system can be replaced locally by another equivalent system which is easier to deal with. A closer examination of this system will lead to an upper bound on the dimension of integral manifolds.
\subsection{Coupled Contact Differential Systems}

It will be more convenient to work in the coordinate system provided by the Lie group coordinates \(Y\in G^-\) of Section~\ref{groupcord}. In this coordinate system the identity matrix corresponds to the reference Hodge structure \(H\). 
Recall that the matrix \(Y\) has a block decomposition \(\sum
Y_{i,j}[i,j]\) where the block \(Y_{i,j}\) at position \((i,j)\) has
size \(h^{w-i,i}\times h^{w-j,j}\). Note that the matrices \(Y_{i,j}\)
must satisfy the Lie group coordinate version of
Equation~\ref{ortcond}, namely
\begin{equation}\label{grouportcond}
(log(Y))^tQ+Qlog(Y)=0
\end{equation}
As a first step let us write the generators of the differential ideal of the horizontal differential system in terms of the entries of the blocks \(Y_{i,j}\).
\begin{prop}\label{firstreduction}
The differential ideal of the horizontal system is generated algebraically by the 1-form entries of the matrices
\begin{equation}\label{almostcontact}
dY_{i,j}-Y_{i,j+1}dY_{j+1,j}\;\;\;\;\rm{for}\;\;\;\;i\geq j+2
\end{equation}
\end{prop}
\begin{proof}
By Proposition~\ref{maurer} we have to show that the entries of the matrices \(\Omega_{i,j}\) for \(i\geq j+2\) can be generated by the 1-forms in Equation~\ref{almostcontact}. We will proceed by induction on \(i-j\). If \(i-j=2\) then an easy computation shows that 
\begin{equation}\label{contactfirst}
\Omega_{i,j}=(Y^{-1}dY)_{i,j}=dY_{i,j}-Y_{i,j+1}dY_{j+1,j}
\end{equation}
so the result holds in this case. Let us now assume that the entries of \(\Omega_{i,j}\) are generated by the forms in Equation~\ref{almostcontact} for \(i-j<k\) and consider \(\Omega_{a,b}\) with \(a-b=k\). Let us compute the \((a,b)\) blocks of both sides of the the defining equation of the Maurer-Cartan matrix 
\[dY=Y\cdot\Omega .\]
We get 
\[dY_{a,b}=\sum_{l=0}^{a}Y_{a,l}\Omega_{l,b}=Y_{a,b+1}dY_{b+1,b}+\sum_{l=b+2}^{a-1}Y_{a,l}\Omega_{l,b}+\Omega_{a,b}.\]
So 
\[\Omega_{a,b}=dY_{a,b}-Y_{a,b+1}dY_{b+1,b}-\sum_{l=b+2}^{a-1}Y_{a,l}\Omega_{l,b}\]
and since in the sum we only have \(\Omega_{l,b}\) with \(l-b\leq a-1-b=k-1\) the right hand side of the above equation is generated by forms in Equation~\ref{almostcontact} by the induction hypothesis and this implies the result.
\end{proof}

As the next step let us define a new differential system on a submanifold of the period domain in a neighborhood of the reference Hodge structure. Let \(W\subset D\) denote the submanifold defined by the equations
\begin{equation}
Y_{i,j}=0 \;\;\;\;\rm{for}\;\;\;\;i>j+2.
\end{equation}

\begin{defi}\label{contactdef}
The coupled contact system on the submanifold \(W\) is given by the differential 1-form entries of the matrices
\begin{equation}\label{contactsystem}
(Y^{-1}d\,Y)_{i,j}=dY_{i,j}-Y_{i,j+1}dY_{j+1,j}\;\;\;\;\rm{for}\;\;\;\;i=j+2. 
\end{equation}
\end{defi}
\begin{rem}
Let us emphasize that the blocks \(Y_{i,j}\) in the previous definition are of course subject to the condition specified by Equation~\ref{grouportcond}.
\end{rem}
\begin{rem}\label{hhell}
Note that the matrices in the definition are exactly the blocks \(\Omega_{i,j}\) for \(i=j+2\) of the Maurer-Cartan matrix written in terms of the Lie group coordinates, i.e., the blocks which are two steps below the diagonal. In other words we defined a differential systems which is the projection of the horizontal system to the submanifold \(W\) (cf. Equation~\ref{almostcontact}). What is remarkable is that this system is locally equivalent to the horizontal system as we will see shortly.
\end{rem}
\begin{rem}
Each of the equations in Definition~\ref{contactdef} is a matrix
valued contact system. In fact, in the weight two case
when \(h^{2,0}=2\) this system is the classical contact system. Note
furthermore that each matrix \(Y_{i,j}\) appears in two
consecutive matrix valued contact systems and so these systems are
coupled through these matrices, hence the name in the definition.
\end{rem}
\begin{thm}\label{contactequiv}
There is a one--to--one dimension preserving correspondence between germs of integral manifolds of the horizontal differential system on \(D\) and germs of integral manifolds of the coupled contact system on \(W\). The correspondence also identifies integral elements of the two systems.
\end{thm}
\begin{rem}
The theorem implies that instead of studying the uniqueness and dimension properties of germs of variations of Hodge structure we can study the corresponding properties of integral manifolds of the coupled contact system. The results we arrive at will remain true for variations of Hodge structure.
\end{rem}
\begin{proof}
By Proposition~\ref{firstreduction} the horizontal system is given by the differential form entries of the matrices
\begin{equation}\label{variation}
dY_{i,j}-Y_{i,j+1}dY_{j+1,j}\;\;\;\;\rm{for}\;\;\;\; i\geq j+2
\end{equation}
and by Definition~\ref{contactdef} the coupled contact system is given by
\begin{equation}\label{contact}
dY_{i,j}-Y_{i,j+1}dY_{j+1,j}\;\;\;\;\rm{for}\;\;\;\; i=j+2.
\end{equation}
Let \(S\subset D\) be a germ of an integral manifold of the horizontal system. Then because of the above equations the projection \(pr_{W}(S)\) to the submanifold \(W\) defines an integral manifold of the contact system. To prove the theorem we have to show that we can go the other way, namely given an integral manifold \(T\subset W\) of the contact system we need to find a {\it locally unique} extension of \(T\) to an integral manifold \(\bar{T}\subset D\) of the horizontal system. This amounts to showing that the differential equations 
\begin{equation}\label{needtosolve}
dY_{i,j}-Y_{i,j+1}dY_{j+1,j}=0\;\;\;\;\rm{for}\;\;\;\; i>j+2
\end{equation}
have unique local solutions, provided that we have a solution of 
\begin{equation}\label{given}
dY_{i,j}-Y_{i,j+1}dY_{j+1,j}=0\;\;\;\;\rm{for}\;\;\;\; i=j+2.
\end{equation}
Since \(T\) is an integral manifold we also have the equations 
\begin{equation}\label{givenwedge}
dY_{i,j+1}\wedge dY_{j+1,j}=0\;\;\;\;\rm{for}\;\;\;\; i=j+2.
\end{equation}
Let us proceed by induction on \(i-j\). Consider the equation 
\begin{equation}\label{solve1}
dY_{a,b}=Y_{a,b+1}dY_{b+1,b}\;\;\;\;\rm{where}\;\;\;\; a-b=3.
\end{equation}
Note that the right hand side of Equation~\ref{solve1} consists of matrices which are already determined at this point since they appear in Equation~\ref{given}. This implies that if we can solve \ref{solve1} for \(Y_{a,b}\) then the solution will be unique. By the Poincar\'e Lemma there is a local solution of \ref{solve1} if and only if the 1-form entries of the right hand side are {\it closed} 1-forms, i.e., it needs to be checked that 
\begin{equation}
dY_{a,b+1}\wedge dY_{b+1,b}=0\;\;\;\;\rm{for}\;\;\;\; a-b=3
\end{equation}
holds. Now, since \(a-(b+1)=2\) we have \(dY_{a,b+1}=Y_{a,b+2}dY_{b+2,b+1}\) by Equation~\ref{given} so
\[
dY_{a,b+1}\wedge dY_{b+1,b}=(Y_{a,b+2}dY_{b+2,b+1})\wedge dY_{b+1,b}=Y_{a,b+2}(dY_{b+2,b+1}\wedge dY_{b+1,b})=0
\]
where the last equality holds by Equation~\ref{givenwedge}. This completes the first step of the induction. 

Assume now that the equations
\begin{equation}\label{givengeneral}
dY_{i,j}=Y_{i,j+1}dY_{j+1,j}
\end{equation}
are solvable for \(i-j<k\). This also implies that 
\begin{equation}\label{wedgegeneralgiven}
dY_{i,j+1}\wedge dY_{j+1,j}=0\;\;\;\;\rm{for}\;\;\;\;i-j<k.
\end{equation}
By the same argument as above we see that to solve 
\begin{equation}\label{solve2}
dY_{a,b}=Y_{a,b+1}dY_{b+1,b}\;\;\;\;\rm{where}\;\;\;\; a-b=k
\end{equation}
we need that
\begin{equation}\label{wedgegeneral}
dY_{a,b+1}\wedge dY_{b+1,b}=0\;\;\;\;\rm{for}\;\;\;\; a-b=k.
\end{equation}
Note again that the matrices on the right hand side of \ref{solve2} have already appeared in Equation~\ref{givengeneral} since \(a-(b+1)=k-1<k\) so uniqueness holds. To check \ref{wedgegeneral}, note that \(dY_{a,b+1}=Y_{a,b+2}dY_{b+2,b+1}\) by \ref{givengeneral} since \(a-(b+1)<k\) so
\[
dY_{a,b+1}\wedge dY_{b+1,b}=Y_{a,b+2}(dY_{b+2,b+1}\wedge dY_{b+1,b})=0
\]
by Equation~\ref{wedgegeneralgiven}.
\end{proof} 
%%%NOTE TO MYSELF: contact is a special case
\subsection{Integral Elements and Abelian Lie Algebras}\label{proofcom}

In this section we will prove Proposition~\ref{commutative}. Let us recall the statement.
\begin{prop}
Let \(E\subset \mathfrak{g}^{-1,1}\) be a linear subspace. Then \(E\) is an integral element of the horizontal differential system if and only if it is an abelian subalgebra.
\end{prop}
\begin{proof}
By definition \(E\) is an integral element if and only if \(\omega_E = 0\) for each \(\omega\) in the differential ideal. By Theorem~\ref{contactequiv} this is equivalent to the vanishing of the two forms \(\varphi_{i,j}=dY_{i,j+1}\wedge dY_{j+1,j}\) on \(E\) for \(i=j+2\). Now let \(X^1,X^2\in E\) be two general tangent vectors. Then 
\begin{equation}
0=\varphi_{i,j}(X^1,X^2)=dY_{i,j+1}\wedge dY_{j+1,j}(X^1,X^2)=X^1_{i,j+1}X^2_{j+1,j}-X^2_{i,j+1}X^1_{j+1,j}
\end{equation}
which is exactly the condition for the commutator \([X^1,X^2]\) to vanish since the matrices \(X^1,X^2\) have nonzero blocks only at the \(i=j+1\) positions.
\end{proof}
\begin{rem}
In the language of Hodge theory integral elements correspond to infinitesimal variations of Hodge structure and the above commutativity condition is part of the definition of an infinitesimal variation of Hodge structure (see \cite{C1}). In this context Proposition~\ref{exptrick} says that to every infinitesimal variation of Hodge structure there is a germ of a variation of Hodge structure tangent to it. Apparently this simple result was not well known.
\end{rem}
\subsection{Dimension Bound for Variations of Hodge Structure}\label{dimbound}

In this section we will use the contact differential system to give sharp upper bounds for the dimension of variations of Hodge structure. This is the main result of \cite{C-K-T}. The proof is different, however, and less elementary than the one in this section. Also a case is not covered in \cite{C-K-T} so it seems worthwhile to reprove this result here.

Up to this point we did not need to use the conditions that Equation~\ref{grouportcond} specifies in terms of the blocks \(Y_{i,j}\) but to give the upper bound we will have to compute these at least for the contact system. By doing this it is possible to further reduce the contact system to an equivalent differential system.
\begin{prop}\label{realreduced}
a) If the weight \(w=2k+1\) is odd then there is a one-to-one correspondence between germs of integral manifolds of the coupled contact differential system and the system given by the 1-form entries of the matrices
\begin{equation}\label{reducedcontactodd}
dY_{i,j}-Y_{i,j+1}dY_{j+1,j}\;\;\;\;\rm{for}\;\;\;\;i=j+2\;\;\;\;\rm{with}\;\;\;\;j<k.
\end{equation}
where the matrix \(Y_{k+1,k}\) is symmetric.

b) If \(w=2k\) is even then the above correspondence is between the coupled contact system and the system given by the entries of
\begin{equation}\label{reducedcontacteven}
dY_{i,j}-Y_{i,j+1}dY_{j+1,j}\;\;\;\;\rm{for}\;\;\;\;i=j+2\;\;\;\;\rm{with}\;\;\;\;j<k-1 \;\;\;\;\rm{and}
\end{equation}
\begin{equation}\label{weight2add}
dY_{k+1,k-1}-Y_{k,k-1}^tdY_{k,k-1}
\end{equation}
and we also have
\begin{equation}\label{weight2cond}
Y_{k+1,k-1}+Y_{k+1,k-1}^t=Y_{k,k-1}^tY_{k,k-1}.
\end{equation}
\end{prop}
\begin{rem}
The content of this proposition is that we need to consider only blocks of the decomposition that are above the main antidiagonal and we have an explicit description of the dependencies between the entries of the matrices (e.g., the condition that \(Y_{k+1,k}\) is symmetric when the weight is odd). To make this easier to visualize let us consider the \(w=3\) and \(w=4\) cases.
\end{rem}
\begin{ex}\label{visex}(\(w=3\))
In this case a general element \(Y\in G^-\) has the form
\[Y=
\begin{pmatrix}
	I & 0 & 0 & 0\\
	Y_{1,0} & I & 0 & 0\\
	Y_{2,0}& Y_{2,1} & I &0\\
	Y_{3,0}& Y_{3,1} & Y_{3,2} &I
\end{pmatrix}
\]
where the blocks satisfy Equation~\ref{grouportcond}. Now the proposition says that the contact differential system is determined by the 1-forms in
\[ dY_{2,0}-Y_{2,1}dY_{1,0}\]
where \(Y_{2,1}\) is symmetric and this is the only dependence between the entries of the matrices giving local coordinates. So to exhibit local variations of Hodge structure in the weight three case it is enough to solve the differential equations in \(dY_{2,0}=Y_{2,1}dY_{1,0}\). 
\end{ex}
\begin{ex}(\(w=4\))
Now a general element \(Y\in G^-\) has the form
\[Y=
\begin{pmatrix} 
	I & 0 & 0 & 0 & 0\\
	Y_{1,0} & I & 0 & 0 & 0\\
	Y_{2,0}& Y_{2,1} & I &0 & 0\\
	Y_{3,0}& Y_{3,1} & Y_{3,2} &I & 0\\
	Y_{4,0}& Y_{4,1} & Y_{4,2} & Y_{4,3}& I
\end{pmatrix}
\]
where the blocks satisfy Equation~\ref{grouportcond}. By the proposition this means that the contact system is determined by
\[ dY_{2,0}-Y_{2,1}dY_{1,0}\;\;\;\;\rm{and}\;\;\;\;dY_{3,1}-Y_{2,1}^tdY_{2,1}\]
where there is one more restriction of the form
\[ Y_{3,1}^t + Y_{3,1}=Y_{2,1}^tY_{2,1}\]
\end{ex}

\begin{proof}
Since the contact system involves blocks that are one or two steps below the main diagonal, to compute these blocks of \(log(Y)\) only the first two terms of the power series expansion will have to be considered: \((Y-I)-(Y-I)^2/2\).
If we consider the blocks that are one step below the diagonal it is immediate from Equation~\ref{ortcond} that
\begin{equation}\label{symmetry}
Y_{i,i-1}^t=Y_{w-i+1,w-i}\;\;\;\;i\in \{0,\ldots ,w\}.
\end{equation}
If \(w\) is odd this implies that \(Y_{k+1,k}\) is symmetric.

Now we have to show that an equation consisting of blocks below the antidiagonal is automatically solved if we have solutions for the equations above the antidiagonal. This is a simple consequence of Equation~\ref{ortcond}. In detail, let us consider the equation 
\begin{equation}
dY_{j+2,j}=Y_{j+2,j+1}dY_{j+1,j}\;\;\;\;\rm{for}\;\;\;\;j>k\;\;\;\;w=2k+1\;\;\;\;\rm{or}\;\;\;\;j>k-1\;\;\;\;w=2k.
\end{equation}
Computing \(log(Y)_{j+2,j}\) and substituting it to Equation~\ref{ortcond} we get
\begin{equation}\label{reducedcontact}
Y_{j+2,j} = Y_{j+2,j+1}Y_{j+1,j}-Y_{w-j,w-j-2}^t
\end{equation}
from which 
\begin{equation}
dY_{j+2,j}=Y_{j+2,j+1}dY_{j+1,j}+dY_{j+2,j+1}Y_{j+1,j}-Y_{w-j,w-j-2}^t=Y_{j+2,j+1}dY_{j+1,j}
\end{equation}
by Equation~\ref{symmetry} and Equation~\ref{reducedcontactodd} (or Equation~\ref{reducedcontacteven}).
\end{proof}

Using this reduced form of the contact system we will proceed to give upper bounds for the dimension of variations of Hodge structure. Let us define the following quadratic functions that depend on the Hodge numbers \(h^{i,j}\). Let \(h^l\) stand for \(h^{l,w-l}\). If the weight is odd \((w=2k+1>1)\)  let
\begin{eqnarray*}
q^{odd}_1& = &\sum_{i=0}^{\infty}h^{k+2+2i}h^{k+3+2i}\cr
q^{odd}_2& = &\frac{1}{2}h^{k+1}(h^{k+1}+1)+\sum_{i=0}^{\infty}h^{k+3+2i}h^{k+4+2i}.
\end{eqnarray*}
If the weight is even \((w=2k)\) let
\begin{eqnarray*}
q^{even}_1& =&\sum_{i=0}^{\infty}h^{k+1+2i}h^{k+2+2i}\cr
q^{even}_2& = &\overline{q_2}+\sum_{i=0}^{\infty}h^{k+2+2i}h^{k+3+2i}\cr
q^{even}_3& = &h^k+h^{k+2}(h^{k+1}-1)+\sum_{i=0}^{\infty}h^{k+3+2i}h^{k+4+2i}\;\;\;\;\rm{if}\;\;\;\;w\geq 4.
\end{eqnarray*}
where 
\begin{equation}
\overline{q_2}=
\left\{\begin{array}{ll}
	h^k & \mbox{if \(h^{k+1}=1\),}\cr
	\frac{1}{2}h^{k+1}h^k& \mbox{if \(h^k\) is even and  \(h^{k+1}>1\),}\cr
	\frac{1}{2}h^{k+1}(h^k-1)+1& \mbox{if \(h^k\) is odd and \(h^{k+1}>1\).} \end{array} \right.
\end{equation}
\begin{rem}
The sums are of course finite since \(h^i=0\) if \(i>w\). \(q^{even}_3\) is defined only if the weight is at least four.
\end{rem}

\begin{thm}\label{dimboundthm}
Let \(S\subset D\) be a variation of Hodge structure.
\begin{enumerate}
\item 
If \(w=2k+1\) then \(dim(S)\leq max\{q^{odd}_1,q^{odd}_2\}\).

\item 
If \(w=2k\) then \(dim(S)\leq max\{q^{even}_1,q^{even}_2,q^{even}_3\}\).
\end{enumerate} 
\end{thm}
\begin{proof}1) According to Proposition~\ref{realreduced} the contact system is equivalent to the system defined by 
\begin{equation}\label{dimproofcontact} dY_{i,j}-Y_{i,j+1}dY_{j+1,j}\;\;\;\;\rm{for}\;\;\;\;i=j+2\;\;\;\;\rm{with}\;\;\;\;j<k.
\end{equation}
with \(Y_{k+1,k}\) symmetric. We will analyze this system to derive the dimension bound. 

In a neighborhood of the reference Hodge structure \(S\) is defined by holomorphic functions depending on a set \(L\) of independent coordinate variables. Clearly \(dim(S)\leq |L|\) where \(|L|\) denotes the cardinality of the set \(L\). We will give a bound for \(|L|\). At the origin the system is defined by the vanishing of the differential form entries of \(dY_{j+2,j}\) \((j<k)\) so the coordinate entries of the matrices \(Y_{j+2,j}\) are not elements of \(L\), i.e.,  all elements of \(L\) appear in the matrices \(dY_{j+1,j}\). Let \(y_{1,0}\) denote the maximal number of the coordinates in the set \(L\) that appear in a column of the matrix \(Y_{1,0}\). Then, clearly, \(Y_{1,0}\) contains at most \(y_{1,0}\cdot h^{w}\) elements of \(L\) since \(Y_{1,0}\) has \(h^{w}\) columns.

Notice that a single entry of the matrix valued differential form \(dY_{2,0}-Y_{2,1}dY_{1,0}\) defines a classical contact system so we can use Proposition~\ref{contactmax} and Remark~\ref{fixj} to partition the variables of this classical contact system into two disjoint subsets. More precisely, suppose that the \(l^{th}\) column in the matrix \(Y_{1,0}\) contains \(y_{1,0}\) elements of \(L\), i.e., it is one of the columns containing the maximal possible number of elements of \(L\). Consider the classical contact system defined by the \((i,l)\) entry of the system \(dY_{2,0}-Y_{2,1}dY_{1,0}\) for every \(i\in \{0,\ldots,h^{w-2}\}\):
\begin{equation}\label{classcont}
(dY_{2,0})^{i,l}-\sum_{t}(Y_{2,1})^{i,t}(dY_{1,0})^{t,l}
\end{equation}
where superscripts denote matrix entries. By
Proposition~\ref{contactmax} the \(i^{th}\) row of the matrix
\(Y_{2,1}\) can contain at most \(h^{w-1}-y_{1,0}\) elements of
\(L\). This is because the classical contact system (\ref{classcont}) already contains
\(y_{1,0}\) independent variables coming from the \(l^{th}\) column of
\(Y_{1,0}\), (the set \(J\) in Proposition~\ref{contactmax}), so since
\(I\) and \(J\) are disjoint we get \(|I|\leq
h^{w-1}-y_{1,0}\). This implies that there are at most
\(h^{w-1}-y_{1,0}\) columns of the matrix \(Y_{2,1}\) that contain
elements of \(L\) so if \(y_{2,1}\) denotes the maximal number of the
coordinates in the set \(L\) that appear in a column of the matrix
\(Y_{2,1}\) then this matrix can contain at most
\(y_{2,1}\cdot (h^{w-1}-y_{1,0}) \) elements of \(L\). 

Applying the same argument to all of the matrix valued contact systems
in Equation~\ref{dimproofcontact} and taking into consideration that
\(Y_{k+1,k}\) is symmetric, we arrive at the following upper bound for \(|L|\):
\begin{eqnarray*}\label{boundd}
|L|& \leq & y_{1,0}\cdot h^{w}+y_{2,1}\cdot
(h^{w-1}-y_{1,0})+\ldots+\\
& + & y_{k,k-1}\cdot (h^{w-(k-1)}-y_{k-1,k-2})+\frac{1}{2}(h^{w-k}-y_{k,k-1})(h^{w-k}-y_{k,k-1}+1)
\end{eqnarray*}

The next step is to determine the maximum value of the right hand side
of the above equation for \(y_{j+1,j}\in [1,h^{w-(j+1)}]\)
(\(j<k)\). This is a quadratic programming problem in a given
rectangle that can be solved as follows. Let \(f(y_{i,j})\) denote
the right hand side of this equation. Then the function \(f(y_{i,j})\) can not have interior maximum because its Hessian matrix is not negative definite and it is nowhere negative semidefinite. If we restrict \(f\) to a face than we get a function which is linear or has a non-negative Hessian. Applying this repeatedly we can see that the maximum must occur at a vertex and examining the function \(f(vertex)\) we conclude the theorem in the odd weight case.

2) First let us consider the \(w=2\) case. By Proposition~\ref{realreduced} the contact system is defined by the equation
\begin{equation}\label{w2}
dY_{2,0}=Y_{1,0}^tdY_{1,0}
\end{equation}
which is subject to the condition
\begin{equation}\label{weight2cond2}
Y_{2,0}+Y_{2,0}^t=Y_{1,0}^tY_{1,0}.
\end{equation}
 \(Y_{1,0}\) determines the symmetric part of \(Y_{2,0}\) by Equation~\ref{weight2cond2} so it is enough to consider the 1-form entries \((i,j)\) of Equation~\ref{w2} for which \(i<j\). The rest of the equations are automatically satisfied as it can be seen by taking the exterior derivative of Equation~\ref{weight2cond2}. 

Let the set \(L\) be as in the first part of the proof and let \(y_{1,0}\) denote the maximal number of the coordinates in the set \(L\) that appear in a column of the matrix \(Y_{1,0}\). (Just like above we can assume that the entries of the matrix \(Y_{2,0}\) do not appear in \(L\)). Then we have the following bounds for \(L\):
\begin{equation}
|L|\leq y_{1,0}\cdot h^{2,0}
\end{equation}
since the matrix \(Y_{1,0}\) has \(h^{2,0}\) columns. We also have 
\begin{equation}
|L|\leq y_{1,0}+(h^{1,1}-y_{1,0})\cdot (h^{2,0}-1)
\end{equation}
which follows by fixing a column having \(y_{1,0}\) coordinate entries in \(L\) and applying the argument of Proposition~\ref{contactmax} to each of the classical contact systems determined by this fixed column and the other columns of \(Y_{1,0}\). Then in the fixed column we have \(y_{1,0}\) independent variables and each of the other columns can have at most \(h^{1,1}-y_{1,0}\) of them. The above two bounds for \(L\) imply that the maximum occurs if \(y_{1,0}=\frac{1}{2}h^{1,1}\). If \(h^{1,1}\) is even then this implies that \(|L|\leq \frac{1}{2}h^{1,1}\cdot h^{2,0}\). If \(h^{1,1}\) is odd then an easy computation shows that \(\frac{1}{2}h^{2,0}(h^{1,1}-1)+1\) is the largest value satisfying both of the bounds. This proves the theorem in the weight two case.

Assume now that \(w\geq 4\). The argument will be the same as in the odd weight case except that when considering the last of the matrix valued differential forms in Equation~\ref{reducedcontacteven} it will be necessary to use similar arguments as we did in the weight two case. Let us consider the last two equations:
\begin{equation}\label{4one}
dY_{k,k-2}=Y_{k,k-1}dY_{k-1,k-2}
\end{equation}
and
\begin{equation}\label{4two}
dY_{k+1,k-1}=Y_{k,k-1}^tdY_{k,k-1}.
\end{equation}
As in the weight two case we need to consider only the antisymmetric part of Equation~\ref{4two}. Using the same arguments as above we have the following bounds for the number of independent variables in the matrix \( Y_{k,k-1}\), (\(dim(Y_{k,k-1})\) will denote this number):
\begin{equation}
dim(Y_{k,k-1})\leq y_{k,k-1}\cdot (h^{k+1}-y_{k-1,k-2}).
\end{equation}
which follows from the fact that by Equation~\ref{4one} \( Y_{k,k-1}\) can have at most \(h^{k+1}-y_{k-1,k-2}\) columns that contain elements of \(L\)  and by definition one column can contain at most \(y_{k,k-1}\) independent variables. We also have 
\begin{equation}
dim(Y_{k,k-1})\leq y_{k,k-1} + (h^{k+1}-y_{k-1,k-2}-1)(h^{k}-y_{k,k-1})
\end{equation} 
from Equation~\ref{4two} by applying the argument we had in the weight two case. 
Now we have to distinguish two cases.\newline
a) If \(h^{k+1}-y_{k-1,k-2}-1\neq 0\) then from the above two bounds we get that the maximum occurs at \(y_{k,k-1}=\frac{1}{2}h^{k}\) which means that 
\begin{equation}
dim(Y_{k,k-1})\leq \frac{1}{2}h^{k}(h^{k+1}-y_{k-1,k-2})\;\;\;\;\rm{if}\;\;\;\;h^{k}\;\;\;\;\rm{even}
\end{equation}
\begin{equation}
dim(Y_{k,k-1})\leq \frac{1}{2}(h^{k}-1)(h^{k+1}-y_{k-1,k-2})+1\;\;\;\;\rm{if}\;\;\;\;h^{k}\;\;\;\;\rm{odd}.
\end{equation} 
Applying arguments as in the odd weight case to the remaining matrix valued contact systems we arrive at the following upper bound for \(L\):

\begin{eqnarray*}
|L| & \leq & y_{1,0}\cdot h^{w}+y_{2,1}\cdot
(h^{w-1}-y_{1,0})+\ldots+\\
& + & y_{k-1,k-2}\cdot (h^{w-(k-2)}-y_{k-2,k-3})+dim(Y_{k,k-1})
\end{eqnarray*}
Now the usual quadratic programming argument applies to this function to give the result in this case.\newline
b) If \(h^{k+1}-y_{k-1,k-2}-1= 0\) then 
\begin{equation}
dim(Y_{k,k-1})\leq h^{k}
\end{equation}
since there can only be one column in \(Y_{k,k-1}\) that contains elements of \(L\). In this case we get the following upper bound:
\begin{eqnarray*}
|L| & \leq & y_{1,0}\cdot h^{w}+y_{2,1}\cdot
(h^{w-1}-y_{1,0})+\ldots+\\
& + & (h^{k+1}-1)\cdot (h^{w-(k-2)}-y_{k,k-1})+h^{k}
\end{eqnarray*}
and applying quadratic programming again we get the claim of the theorem.
\end{proof} 
\begin{rem}
Let \(\Phi : M \longrightarrow \Gamma~\backslash~ D\) be a global variation of Hodge structure. Then \(rank(\Phi)\) also satisfies the dimension bound given in the theorem.
\end{rem}
\begin{rem}
For sharpness of this result see \cite{C-K-T}. Note that the case when the maximal dimension is \(q^{even}_3\) is missing from the main theorem in \cite{C-K-T}. It is easy to see that there are Hodge numbers for which \(q^{even}_3\) is indeed the maximal dimension so this function is needed in the proper form of the dimension bound.
\end{rem}
%\begin{ex}
%weight 3 and 4
%\end{ex}

\section{Rigidity of Maximal Dimensional  Variations of 
Hodge Structure}

In this section we will examine maximal dimensional variations. We will proceed in two steps. First, we fix a maximal dimensional integral element of the horizontal system and investigate the nature of germs of integral manifolds tangent to this integral element. It turns out that in most cases there is a unique germ of an integral manifold tangent to the given integral element. This is sometimes referred to as local rigidity. In the second step we will fix a point and investigate the nature of integral elements through the given point hence obtaining infinitesimal rigidity. These two steps together lead to the main rigidity results.

\subsection{Germs of Integral Manifolds}

Let \(E\) be a maximal dimensional integral element of the horizontal differential system going through the point \(p\). By Theorem~\ref{exptrick} there is a maximal dimensional integral manifold tangent to \(E\) so the only question we have to investigate is uniqueness of such integral manifolds. The following result answers this question.
\begin{thm}\label{localrigid}
Let \(E\) be a maximal dimensional integral element of the horizontal differential system. 
\begin{enumerate}
\item 
Let \(w=2k+1\). If \(h^{k,k+1}>2\) and all the other Hodge numbers are greater than one then there is a unique germ of an integral manifold whose tangent space at \(x\) is \(E\).

\item Let \(w=2k\). Assume that one of the following conditions holds:
\begin{enumerate}
\item
\(dim(E)=q^{even}_1\) and all the Hodge numbers are greater than one
\item \(dim(E)=q^{even}_2\), \(h^{k+1}>2\), the other Hodge numbers are greater than one and \(h^k\geq 4\) is even,
\item \(dim(E)=q^{even}_3\) and  \(h^{k+1}>2\)
\end{enumerate}
then there is a unique germ of an integral manifold  whose tangent space at \(x\) is \(E\).

\item If \(dim(E)=q^{even}_2\) and \(h^k\) is odd then there is an infinite dimensional family of germs of maximal dimensional integral manifolds tangent to \(E\). 
\end{enumerate} 
\end{thm}
\begin{proof}
1) Assume first that the maximum is \(q^{odd}_1\) i.e., \(dim(E)=q^{odd}_1\). Consider the first matrix valued differential form of Equation~\ref{reducedcontactodd} 
\begin{equation}\label{egy}
dY_{k+1,k-1}=Y_{k+1,k}dY_{k,k-1}
\end{equation}
and its exterior derivative
\begin{equation}\label{ketto}
0=dY_{k+1,k}\wedge dY_{k,k-1}
\end{equation}
From the proof of Theorem~\ref{dimboundthm} we see that in this case the entries of the matrix \(Y_{k,k-1}\) are independent variables (elements of the set \(L\)). We would like to see that the entries of the matrix \(Y_{k+1,k}\) are at most linear functions of the variables in \(L\) since this implies that they are uniquely determined by fixing the tangent space \(E\). To this end, consider an entry of Equation~\ref{ketto} (e.g., the \((1,1)\) entry):
\begin{equation}\label{harom}
0=\sum_{i}(dY_{k+1,k})^{1,i}\wedge (dY_{k,k-1})^{i,1}.
\end{equation}
This equation involves the first row of \(dY_{k+1,k}\) and the first
column of \(dY_{k,k-1}\). Since all the entries of \(dY_{k,k-1}\) are
independent, Cartan's lemma (\cite{W} Exercise 2.16) implies that 
\begin{equation}
(dY_{k+1,k})^{1,i}=\sum_jA_{i,j}(dY_{k,k-1})^{j,1}\;\;\;\;\rm{for}\;\;\;\;\rm{each}\;\;\;\;i
\end{equation}
for some functions \(A_{i,j}\). Since the matrix \(Y_{k,k-1}\) has at least two columns because of the condition on the Hodge numbers in the theorem, we can do the same for the (1,2) entry of Equation~\ref{ketto}. The same argument shows that
\begin{equation}
(dY_{k+1,k})^{1,i}=\sum_jB_{i,j}(dY_{k,k-1})^{j,2}\;\;\;\;\rm{for}\;\;\;\;\rm{each}\;\;\;\;i
\end{equation}
for some functions \(B_{i,j}\). Since the entries in the first and
second column of the matrix \(Y_{k,k-1}\) are independent from each
other this is  possible only if all the functions \(A_{i,j}\) and \(B_{i,j}\) are zero. This implies that entries in the first row of the matrix \(Y_{k+1,k}\) must be constant functions. The same argument applies to all the rows of this matrix so we conclude that all the entries of the matrix \(Y_{k+1,k}\) are constant functions. 

Consider now the next matrix valued differential form of Equation~\ref{reducedcontactodd} 
\begin{equation}\label{egyy}
dY_{k,k-2}=Y_{k,k-1}dY_{k-1,k-2}.
\end{equation}
Notice that this equation contains the matrix \(Y_{k,k-1}\) which has independent entries so applying the above argument we conclude that the entries of the matrix \(dY_{k-1,k-2}\) are constant functions. Continuing in the same manner we see that the matrices involved in Equation~\ref{reducedcontactodd} are either constant matrices or contain entries that are independent coordinate functions on the integral manifold. This implies the claim of the theorem showing that there is a unique germ of an integral manifold tangent to the fixed integral element \(E\). 

If the maximum is \(q^{odd}_2\)  (i.e., \(dim(E)=q^{odd}_2\)) the above argument applies but we have to be careful about the equation containing the symmetric matrix \(Y_{k+1,k}\). Consider this equation 
\begin{equation}\label{hello}
dY_{k+1,k-1}=Y_{k+1,k}dY_{k,k-1}
\end{equation}
and its exterior derivative
\begin{equation}
0=dY_{k+1,k}\wedge dY_{k,k-1}.
\end{equation}
From Theorem~\ref{dimboundthm} we see that the entries of the symmetric matrix \(Y_{k+1,k}\) are independent variables and we would like to conclude that this implies that \(Y_{k,k-1}\) has constant entries. Consider the first column of \(Y_{k,k-1}\). As before from the \((1,1)\) entry of Equation~\ref{hello} we have
\begin{equation}
(dY_{k,k-1})^{i,1}=\sum_jA_{i,j}(dY_{k+1,k})^{1,j}\;\;\;\;\rm{for}\;\;\;\;\rm{each}\;\;\;\;i
\end{equation}
and from the \((2,1)\) entry we have 
\begin{equation}
(dY_{k,k-1})^{i,1}=\sum_jB_{i,j}(dY_{k+1,k})^{2,j}\;\;\;\;\rm{for}\;\;\;\;\rm{each}\;\;\;\;i.
\end{equation}
Since \(Y_{k+1,k}\) is symmetric this  implies only that the entries in the first column of \(Y_{k,k-1}\) are functions of the \((1,2)\) entry of \(Y_{k+1,k}\). On the other hand considering the \((1,1)\) and \((1,3)\) entries of Equation~\ref{hello} (\(h^{k,k+1}>2\)) we conclude that these same entries are functions of the \((1,3)\) entry of \(Y_{k+1,k}\). These two facts together imply that the entries must be constant functions. The same argument applies to all the remaining entries and to all of the remaining matrix valued contact equations and this implies the claim of the theorem in this case. 

2)
a) The difference between the odd weight and even weight cases is the appearance of the equation 
\begin{equation}\label{extraeq}
dY_{k+1,k-1}=Y_{k,k-1}^tdY_{k,k-1}.
\end{equation}
However, if the maximum dimension occurs at \(dim(E)=q^{even}_1\) then there are no independent variables among the entries of matrix \(Y_{k,k-1}\); consequently the same argument applies as in the odd weight, \(q^{odd}_1\) case.

b) Assume now that \(dim(E)=q^{even}_2\). According to Theorem~\ref{dimboundthm} this implies that each column of the matrix \(Y_{k,k-1}\) has exactly \(\frac{1}{2}h^k\) independent entries. Without loss of generality we can assume that the first \(\frac{1}{2}h^k\) entries in the first column (\((Y_{k,k-1})^{l,1}\) with \(l\leq \frac{1}{2}h^k\)) are independent. Consider the \((1,2)\) entry in Equation~\ref{extraeq}:
\begin{equation}
(dY_{k+1,k-1})^{1,2}=\sum_l(Y_{k,k-1})^{l,1}(dY_{k,k-1})^{l,2}.
\end{equation}
Applying Proposition~\ref{contactmax} to this classical contact system it follows that the independent entries in the second column are 
\begin{equation}
(Y_{k,k-1})^{l,2}\;\;\;\;\rm{with}\;\;\;\;l>\frac{1}{2}h^k.
\end{equation}
Furthermore, by Cartan's lemma, the entries \((dY_{k,k-1})^{l,1}\) with \(l>\frac{1}{2}h^k\) can be expressed in terms of the independent 1-forms \((dY_{k,k-1})^{l,1}\) with \(l\leq \frac{1}{2}h^k\) and \((dY_{k,k-1})^{l,2}\) with \(l>\frac{1}{2}h^k\). Considering the \((1,3)\) entry of Equation~\ref{extraeq} we can similarly conclude that the same entries can be expressed in terms of \((dY_{k,k-1})^{l,1}\) with \(l\leq \frac{1}{2}h^k\) and \((dY_{k,k-1})^{l,3}\) with \(l>\frac{1}{2}h^k\). This implies that these entries depend only on \((dY_{k,k-1})^{l,1}\) with \(l\leq \frac{1}{2}h^k\). Similarly, the entries \((dY_{k,k-1})^{l,2}\) with \(l\leq \frac{1}{2}h^k\) depend only on \((dY_{k,k-1})^{l,2}\) with \(l>\frac{1}{2}h^k\). It follows from this and the \((1,2)\) entry of the exterior derivative of Equation~\ref{extraeq} 
\begin{equation}
0=\sum_l(dY_{k,k-1})^{l,1}\wedge (dY_{k,k-1})^{l,2}
\end{equation}
that the entries \((Y_{k,k-1})^{l,1}\) with \(l>\frac{1}{2}h^k\) can be at most linear functions which is what we wanted to prove.

Consider now the exterior derivative of the next matrix valued contact system in Equation~\ref{reducedcontacteven}
\begin{equation}
0=dY_{k,k-1}\wedge dY_{k-1,k-2}
\end{equation}
By applying changes of coordinates we can assume that the independent variable entries in matrix \(Y_{k,k-1}\) are the first \(\frac{1}{2}h^k\) entries in each column. This implies that at least the first two rows of matrix \(Y_{k,k-1}\) consist of independent variables, which in turn implies that the entries of the matrix \(Y_{k-1,k-2}\) must be constant functions by applying the usual Cartan's lemma type argument. All the remaining matrix valued contact systems in Equation~\ref{reducedcontacteven} can be treated exactly the same way as in the odd weight case so we conclude that the entries in these matrices which are not independent variables must be constants. Together with the above results this concludes the proof for this case.

c) In this case \(w\geq 4\) and \(dim(E)=q^{even}_3\). By Theorem~\ref{dimboundthm} the matrix \(Y_{k,k-1}\) has exactly one column that consists entirely of independent variables and the other columns do not contain independent variables at all. We can assume that the entries of the first column are independent. Considering the entries of Equation~\ref{weight2add} it follows that the remaining entries of \(Y_{k,k-1}\) can be expressed as functions of the variables in the first column. In fact, from the \((1,2)\) entry of the exterior derivative of Equation~\ref{weight2add}
\begin{equation}
0=\sum_l(dY_{k,k-1})^{l,1}\wedge (dY_{k,k-1})^{l,2}
\end{equation}
it follows that the entries \((Y_{k,k-1})^{l,2}\) depend only on the independent variables \((Y_{k,k-1})^{l,1}\). Similarly, the same holds for the other columns.

In this case the matrix \(Y_{k-1,k-2}\) consists of independent variables except for its first row (cf. proof of Theorem~\ref{dimboundthm}). From the \((1,1)\) entry of the exterior derivative of Equation~\ref{4one}  
\begin{equation}
0=\sum_l(dY_{k,k-1})^{1,l}\wedge (dY_{k-1,k-2})^{l,1}
\end{equation}
we can conclude that the entries \((Y_{k,k-1})^{1,l}\) for \(l\geq 2\)
depend only on \((Y_{k,k-1})^{1,1}\) and in fact, they must be linear
functions in this variable. Similarly, the entry
\((Y_{k-1,k-2})^{1,1}\) must be a linear function of the variables
\((Y_{k-1,k-2})^{1,l}\) \((l\geq 2)\). These claims easily follow from
the fact that the differential 1-forms \((dY_{k,k-1})^{1,1}\) and
\((dY_{k-1,k-2})^{1,l}\), \((l\geq 2)\) are linearly
independent. Considering the other entries we conclude that the
entries that are not independent variables must be linear functions. In fact, let us remark here that the coefficients of the linear functions in a given row of the matrix \(Y_{k,k-1}\) must be the same for each entry since for example in the second column each coefficient is equal to \(\frac{\partial((Y_{k-1,k-2})^{1,1})}{\partial((Y_{k-1,k-2})^{2,1})}\). It also follows that the coefficients of the linear functions \((Y_{k-1,k-2})^{1,i}\) are the same for each \(i\).

The remaining matrix valued contact systems in Equation~\ref{reducedcontacteven} can now be treated the same way as in the odd weight case which concludes the proof of the theorem in this case.

3) What remains to be considered to complete the proof of the theorem  is the case when \(h^k\) is odd. We will exhibit an infinite dimensional family of germs of maximal dimensional integral manifolds tangent to an integral element \(E\). To this end let us specify the entries of the matrices in Equation~\ref{reducedcontacteven}. To give a maximal dimensional integral manifold we must specify \(\frac{1}{2}h^{k+1}(h^k-1)+1\) independent variables among the entries of \(Y^{k-1,k}\), according to Theorem~\ref{dimboundthm}. Let us choose the first \(\frac{1}{2}(h^k-1)\) variables in each column of \(Y^{k-1,k}\) to be independent and let the next \(\frac{1}{2}(h^k-1)\) entries in each column to be \(\sqrt{-1}\) times the first \(\frac{1}{2}(h^k-1)\) independent entries. This leaves the last row to be considered. Let the first element of the last row be an independent variable and let the remaining entries of the last row be arbitrary holomorphic functions of of the first entry. 

For example, if \(h^k=5\) and \(h^{k+1}=4 \) then 
\[
Y_{k-1,k}=
\begin{pmatrix}
	y_{1,1}& y_{1,2} & y_{1,3} & y_{1,4}\\
	y_{2,1}& y_{2,2} & y_{2,3} & y_{2,4}\\
	i\cdot y_{1,1}& i\cdot y_{1,2} & i\cdot y_{1,3} & i\cdot y_{1,4}\\
	i\cdot y_{2,1}& i\cdot y_{2,2} & i\cdot y_{2,3} & i\cdot y_{2,4}\\
	y_{5,1}& f_{5,2} & f_{5,3} & f_{5,4}
\end{pmatrix}
\]
where the functions \(f_{a,b}\) are arbitrary holomorphic functions of the single variable \(y_{5,1}\) and \(i=\sqrt{-1}\).

Furthermore, let the matrices \(Y_{k-3-2j,k-2-2j}\) (for each \(j\)) consist entirely of independent variables and let the remaining matrices \(Y_{k-2-2j,k-1-2j}\) (for each \(j\)) be zero. This means that we have \(q^{even}_2\) independent variables. It is an easy computation to verify that all the equations in (\ref{reducedcontacteven}) are satisfied so we have defined a maximal dimensional integral manifold. Fixing an integral element to which this integral manifold has to be tangent to can determine only the linear parts of the functions \(f_{a,b}\) which implies that by varying these functions we can exhibit an infinite dimensional family of germs of integral manifolds tangent to a fixed integral element.  
\end{proof}
\subsection{Maximal Dimensional Integral Elements}

In this section we examine how maximal dimensional integral elements are related to each other. First let us recall a theorem of Carlson (\cite{C1}) that will be used in this section.
\begin{thm}[Carlson]
Let \(D\) be the period domain of weight two Hodge structures with \(h^{2,0}>2\) and \(h^{1,1}\) even. Let \(E_1\) and \(E_2\) be maximal dimensional integral elements of the horizontal distribution. Then there is an element \(g\in Aut(D)\) such that \(g\cdot E_1=E_2\).
\end{thm}
We will consider the same question for higher weight, namely what can be said about the relationship between maximal dimensional integral elements of the horizontal differential system. 
\begin{thm}\label{elementrigid}
Let \(D\) denote the period domain. 
\begin{enumerate} 
\item Assume that one of the following holds:
\begin{enumerate}
\item \(w=2k+1\), \(h^{k,k+1}>2\) and all the other Hodge numbers are greater than one
\item \(w=2k\), all the Hodge numbers are greater than one and the maximum dimension is \(q^{even}_1\) 
\end{enumerate}
then there is a unique maximal dimensional integral element of the horizontal system through each point of \(D\).
\item Let \(E_1\), \(E_2\) be two maximal dimensional integral elements. Assume that one of the following holds:
\begin{enumerate}
\item \(w=2k\), \(dim(E_1)=dim(E_2)=q^{even}_2\), \(h^{k+1}>2\), the other Hodge numbers are greater than one and \(h^k\geq 4\) is even
\item \(w=2k\), \(dim(E_1)=dim(E_2)=q^{even}_3\) and \(h^{k+1}>2\) 
\end{enumerate}
then there is an element \(g\in Aut(D)\) such that \(g\cdot E_1=E_2\).
\end{enumerate}
\end{thm}
\begin{proof}
1) a) and b) The horizontal system is homogeneous so we can always assume that the integral elements are at the reference Hodge structure. In Theorem~\ref{localrigid} we proved that in cases 1.a and 1.b the matrix entries are either independent variables or constants, and these constants are determined be the point of the period domain. If we are at the reference Hodge structure then all these constants are zero. To compute the tangent space to the integral manifold we have to take partial derivatives of the matrix entries by the independent variables. From this it is clear that there can be only one maximal dimensional integral element through each point of \(D\).

2) a) In this case it follows from the proof of Theorem~\ref{localrigid} that the matrices besides \(Y_{k,k-1}\) consist of independent variables or constant functions, so the subspace of the tangent space coming from these matrices is unique. Now, it is clear that the equations \(Y_{k+1,k}\) has  to satisfy are exactly the weight two equations so Carlson's theorem applies to conclude the theorem in this case. 

2) b) Again, from the proof of Theorem~\ref{localrigid} we see that the subspace of the tangent space coming from the matrices besides \(Y_{k,k-1}\) and \(Y_{k-1,k-2}\) is unique, so we need to consider these two matrices. About these matrices we proved that they contain linear entries such that the coefficients of the columns are the same. This implies that we can conjugate the possible tangent spaces into each other by multiplying by elements of \(Aut(D)\).
\end{proof}
\begin{rem}
It is true that result 2) of the above theorem holds even if \(h^{k}\) is odd. However, the author does not know a simple proof of this along the lines of these other results. Since there is an infinite dimensional family of integral manifolds to a given integral element in this case, we could not use this result to conclude further rigidity properties; hence it is not proved here. 
\end{rem}
\begin{rem}
Theorem~\ref{localrigid} and Theorem~\ref{elementrigid} immediately imply our main result Theorem~\ref{mainthm}. Note that if the Hodge numbers are smaller than what is required for these results, then we are reduced to the classical contact case and so flexibility holds by Theorem~\ref{arnold}. 
\end{rem}
\begin{rem}
It is not known whether the maximal dimensional variations are geometric if the weight is bigger than two. In the weight two, \(h^{1,1}\) even case this is true by \cite{C-S}.

It is an interesting question in the author's opinion whether the reduction of the horizontal system to the much smaller coupled contact system is merely a local possibility or there is some global geometric reason that would explain this result.
\end{rem}

%%%%%%%%%%%%%

\end{document}